# A scale-based approach to finding effective dimensionality in manifold learning

### Xiaohui Wang[*]

*Department of Statistics*
*University of Virgina*
*e-mail:* xw5a@virginia.edu

### J. S. Marron[†]

*Department of Statistics and Operation Research*
*University of North Carolina*
*e-mail:* marron@email.unc.edu

**Abstract:**
The discovering of low-dimensional manifolds in high-dimensional data is one of the main goals in manifold learning. We propose a new approach to identify the effective dimension (intrinsic dimension) of low-dimensional manifolds. The scale space viewpoint is the key to our approach enabling us to meet the challenge of noisy data. Our approach finds the effective dimensionality of the data over all scale without any prior knowledge. It has better performance compared with other methods especially in the presence of relatively large noise and is computationally efficient.

## 1. Introduction

High-dimensional data sets can have meaningful low-dimensional structures hidden in the observation space. However, a data set can have different effective dimensions at different scales in the presence of noise. Figure 1 shows two toy data sets. Plot (a) shows the data randomly drawn from the well-known Swiss roll. The Swiss roll is a 2-d submanifold embedded in a 3-d space and can be thought of as curling a piece of rectangular paper. Plot (b) shows the data randomly drawn from the same Swiss roll plus 3-d Gaussian noise. The noisy Swiss roll data can be viewed as either 2-d or 3-d depending on scale. At coarse scales the noise is negligible, so the data are essential 2-d. At fine scales the 3-d noise dominates, so no 2-d structure is present.

---

[*]X. Wang is an assistant professor of Department of Statistics at University of Virginia and this research was supported by the NSF Graduate Student Fellowship and partially supported by the NSF grants DMS-0631639.

[†]J. S. Marron is a professor of Department of Statistics and Operation Research at University of North Carolina at Chapel Hill and this research was partially supported by the NSF grants DMS-9971649 and DMS-0308331. This work is from the Ph.D. dissertation of the first author under the supervision of the second author.

(a) Swiss Roll without noise     (b) Swiss Roll with noise

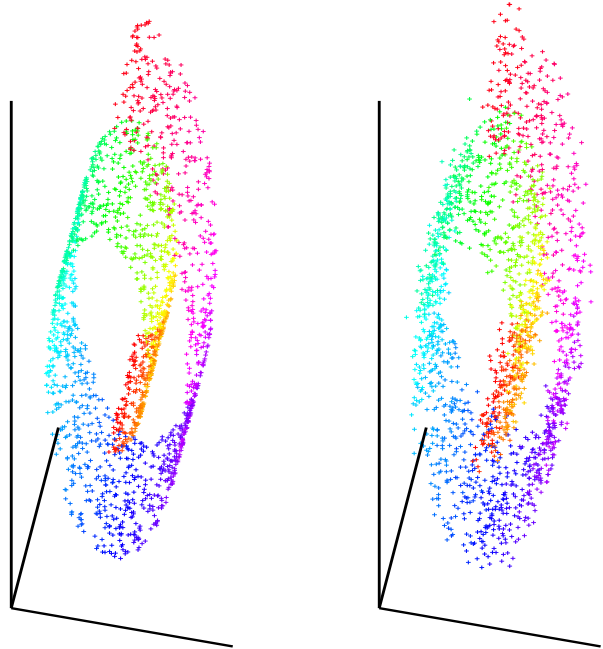

FIG 1. *Examples of 3-d data lying near the 2-d manifold. (a) data on the Swiss roll. (b) data on the same Swiss roll plus 3-d Gaussian noise. The latter is particularly challenging for conventional dimension methods but easily handled by our scale-space approach*

In this paper, we will not endeavor to estimate the true dimension, but to find an *effective dimension* (or intrinsic dimension as in many other literatures), i.e., the dimension that will give a reasonable fit. We realize that the effective dimension can be chosen by many other methods, such as generalized cross-validation (GCV) and cross-validation (CV). However, they are computationally intensive in practice and have some well documented poor realistic issues, see Jones, Marron and Sheather (1996). Here we propose a simple direct estimate which would serve as an input for many methods described later. This is done through a set of hypothesis tests to extract the *effective dimensions* of data. The challenge of noisy data is tackled using a scale space approach. The effective dimensions of the data are estimated over all scales without any prior knowledge, which allows a much larger amount of noise than earlier methods can handle. The multi scale estimated effective dimensionality of the data not only gives information about the lowest dimension of the submanifold on which the data lie closely, but also indicates the noise level present in the data.

In the past decades, Principal Component Analysis (PCA) and related methods have been workhorse approaches to data sets whose low-dimensional underlying structure is linear or can be approximated linearly. Several approaches have been devised to address the problem of finding low-dimensional nonlinear structure in data. Principal Curve Analysis (Hastie and Stuetzle 1989) mainly focused on the 1-dimensional nonlinear structure. Two-dimensional surfaces are investigated by Hastie (1984). LeBlanc and Tibshirani (1994) extended the idea to higher dimension cases. In their paper, they pointed out that an important aspect is the proper lowest dimension which the model was based on. They used the generalized cross-validation with the backward stepwise pruning algorithm to select the dimension.

Some methods of Nonlinear PCA have been also developed. Gnanadesikan (1977) applied PCA to a vector of which components are polynomial terms generated from the data. The method has difficulty with high-dimensional data, since the number of polynomial terms increases rapidly with the dimensionality. Kernel PCA by Schölkopf et al. (1998) is widely developed. However, the result is very hard to interpret and gain meaningful insight from. The Sandglass-type Multi-Layered Perceptron by Irie and Kawato (1990) and by Demers and Cotterell (1993) can construct adequate nonlinear mapping functions to extract a low-dimensional internal representation from given data. However, the dimensionality of the internal representation is fixed and depends on previous knowledge of data. Again, determining the intrinsic dimension is one of the key points.

Recently, a variety of methods have been developed to deal with nonlinear dimensionality reduction. Among them, they are Isometic Feature Mapping (ISOMAP) (Tenenbaum, de Silva, and Langford 2000), Local Linear Embedding (LLE) (Roweis and Saul 2000), Hessian-based Locally Linear Embedding (Donoho and Grimes 2003), and others. Those methods focus on finding a low-dimensional curved manifold embedding of high-dimensional data. The dimensionality of the embedding is a key parameter; however, there is no consensus on how such dimensionality is determined. The dimensionality has been heuristically chosen from the curve of residual variance as a function of dimension. Constructing a reliable estimator of the intrinsic dimension and understanding its statistical properties will clearly improve the performance of manifold learning methods.

The current dimension estimating methods can be roughly divided into two groups, the eigenvalue methods and the geometric methods. Eigenvalue methods are based on either PCA (Fukunaga and Olsen 1971) or local PCA (Bruske and Sommer 1998). PCA can be very ineffective for nonlinear data. For example, applying PCA to the data in Figure 1 will show that a 3-dimension representation is needed to represent them. Local PCA depends heavily on the choice of local region and threshold. When the data with non-linear underlying structure, the choice of the local region and threshold depends on the noise level. Without any prior knowledge, it is hard for local PCA to get a reasonable estimation of true dimensionality of the data. The method developed in this paper serves as an exploratory tool for discovering the dimensionality of the data without any prior knowledge of the noise level. If the underlying structure of the data is linear, our method is consistent with PCA. When the underlying structure is non-linear, our method is more effective than local PCA when the noise level is unknown.

The geometric methods are mostly based on fractal dimensions or nearest neighbor distances. Details can be found in Grassberger and Procaccia (1983); Camastra and Vinciarelli (2002); Costa and Hero (2004). The statistical properties have been studied in Levina and Bickel (2005). Smith (1992) discussed how difficulty it is to estimate the fractal dimension in noisy chaotic time series. In Section 4.3, we give examples to show that our method is more robust than the methods of the standard fractal dimension estimation.

The method developed in this paper has been implemented in Matlab code, which is also available from authors' web[1]. For all the simulated examples and real data examples in this paper, applying our method only takes less than a minute on a simple personal computer. The Matlab program will generate the estimated effective dimensions for all scale, which is straight forward to interpret.

The paper is organized as follows. Since we compare our method mostly with the estimated dimensionality by ISOMAP, we will first introduce ISOMAP in Section 1.1. This is followed in Section 2 by the derivation of our procedure, including the scale space idea. In Section 3, we define *Vector Dimensionality* which is related to the population effective dimensionality, set up a set of hypothesis tests, prove the consistency of the test statistic. Section 4 contains examples. Section 5 gives summary and further discussion.

## *1.1. ISOMAP*

Given the distances (similarities) of pairwise data points, Multidimensional Scaling (MS) techniques try to find a representation of the data in low dimensions such that the distances (inter-item proximities) in low-dimensional representation space nearly match the original distances (similarities). ISOMAP builds on the MS algorithm but uses geodesic distance between all pairs of data points. Figure 1 (a) shows that points far apart on the Swiss roll have short Euclidean distances. The key point of ISOMAP is to capture the geodesic manifold distance between all pairs of data points. ISOMAP estimates the geodesic distances for neighboring points by their Euclidean distances and the geodesic distances for far apart points by finding the shortest path in a graph with edges connecting neighboring data points. In addition to giving a low-dimensional representation, ISOMAP indicated potential estimated dimensionality of the data through an interpretive error curve as a function of dimensionality.

When data points lie exactly on the manifold or have relatively little noise, the neighboring points based on Euclidean distance are consistent with the neighboring points based on the manifold geodesic distance. Balasubramanian and Schwartz (2002) argued the topological instability of ISOMAP when data have relatively large noise, since ISOMAP fails to give a proper estimation of intrinsic dimension for the latter case. In 1974, Shepard addressed a method which had the similar idea as ISOMAP. Shepard also pointed out that it was feasible to find the intrinsic dimension of the data particularly in cases of relatively noise-free data. One important drawback is that ISOMAP fails when the underlying manifold is not isomophic to Euclidean space. For example, if the underlying structure of the data is a circle, ISOMAP will fail to determine its true dimension. This will be illustrated in Example 1 of Section 4.

---

[1] The Matlab code is available at faculty.virginia.edu/xiaohui

## 2. Scale space approach

The notion of scale has been deeply studied in the field of computer vision. For an introduction and detailed discussion, see Lindeberg (1993) and ter Haar Romeny (2002). The problem of scale must be faced in any imaging situation. An inherent property of objects in the world and details in images is that they only exist as meaningful entities over certain ranges of scale. This is done by studying a family of Gaussian window smooths of the image, indexed by the window width. While smoothing of images has been done in many contexts, what distinguishes the scale space approach is the idea of considering all scales, instead of trying to choose a level of smoothing. This approach is very useful because often different aspects of the underlying signal can be most clearly seen at different levels of scale. We view the problem of finding meaningful low-dimensional structure in high dimensional data as analogous to finding meaningful image properties. When data have noise, the meaningful underlying structure depends on the noise level. Studying a range of different scales helps to detect the meaningful low-dimensional structure from the high-dimensional noise. In this section, we introduce our notion of scale parameter. To fully understand the underlying structure of the data, it is quite important to study useful statistics across all scales.

Figure 2 shows three toy data sets which are uniformly generated on a line segment (with length equal to 1) plus vertical Gaussian noise realizations with different standard deviation ($\sigma = 0.001, 0.03, 0.3$ from left to right). The two different circles are two different window sizes. In Figure 2, there are clear differences among these three cases. The left panel shows data that are very close to lying on a line, a 1-d submanifold (i.e., lower-dimensional structure). The right panel shows data that are not close to a 1-d manifold but are much more "2-dimensional". The middle panel is between. All three toy data sets, as shown in Figure 2, show different high-dimensional structures at the small scale with respect to the different noise level. The clarity of the low-dimensional underlying structure (the horizontal line) is affected in different ways by the noise at the different scales. We will introduce a new parameter, scale $s$. Through it, we can detect the effective dimension of the data even in the presence of noise.

In the presence of noise, none of these toy data sets lie exactly on a 1-d manifold. However, it is clear that some are much closer than others to lying on a 1-d manifold. To measure this closeness, we use the scale space idea. Two candidate scales are illustrated by the big and small circles in Figure 2. Viewing the data points through the big window, the data set in the left panel is close to lying on a line, and the one in the middle panel is close to lying on a thick line. Viewing the data points through the small window, we may still claim that the data set in the left panel nearly lies on a line, but the one in the middle panel is more "2-dimensional". For both such window sizes, the data set in the right panel is more like lying on the plane instead of a line. From this example, notice that for a given data set with noise, whether it is close to a low-dimensional submanifold or not depends on the window size.

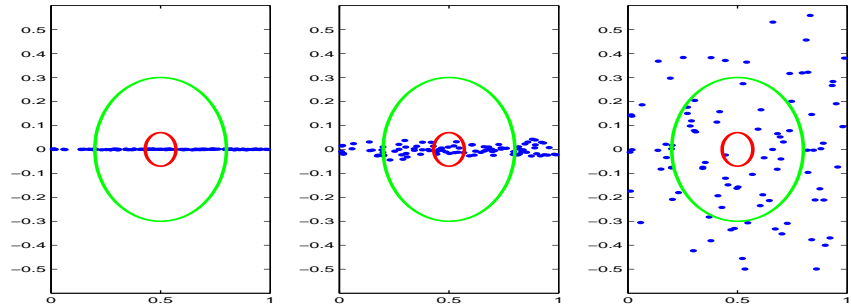

FIG 2. *Three panels show toy data sets with each sample size equal to 100 generated uniformly on a line segment with length equal to 1 plus different independent vertical noise realizations. The co-centered circles represent the different scales. The 3 examples provide the different challenges to the effective different dimensions at the scales.*

To deal with different window sizes, a new scale parameter, $s$, is introduced, which is the radius of the circle for the toy examples in Figure 2. If the data are from $d$-space, the scale parameter $s$ is defined as the radius of the $d$-dimensional ball. The scale parameter $s$ essentially measures the window size. Here the radius of the circle is just one choice. Since the underlying structure of the noisy data changes as the window size changes, any information used to detect the low dimensional structure will be a function of $s$ for $0 < s < \infty$. The lesson we learned from the examples in Figure 2 is: (i) data exhibit a high-dimensional structure when the scale $s$ is less than the noise level; (ii) when the scale $s$ is greater than the noise level, it is possible to find meaningful low dimensional structure. However, a key issue is that the noise level is in general unknown. Without any prior information about the data, useful insights can be obtained from considering a range of scales. Therefore, we recommend studying the data at all scales to get the whole picture.

### *2.1. Dimension test statistic*

How can we characterize a data set that lies on a 1-d manifold? For any distinct three points on the straight line (a 1-d linear submanifold) $\mathbf{x}_i$, $\mathbf{x}_j$ and $\mathbf{x}_k$, the angle between the vector $\mathbf{x}_j - \mathbf{x}_i$ and the vector $\mathbf{x}_k - \mathbf{x}_i$ is either 0 or $\pi$. This is in general not true when the 1-d submanifold is curved. However, if we consider the two points which are near to $\mathbf{x}_i$ for data lying on a curved 1-d submanifold, the corresponding angle should also be close to either 0 or $\pi$.

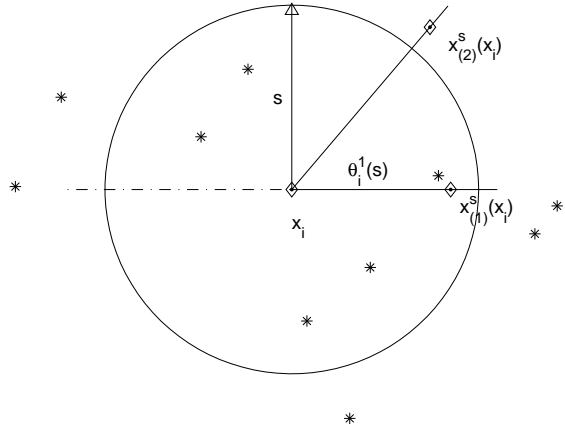

FIG 3. *Toy data set illustrating* $\mathbf{x_i}$, $\mathbf{x^s_{(1)}}(\mathbf{x_i})$, $\mathbf{x^s_{(2)}}(\mathbf{x_i})$, *and how the angle* $\theta^1_i(s)$ *is formed with data points in 2-d space.*

For each point $\mathbf{x}_i$, and a fixed value of scale $s$ with $0 < s < \infty$, first define $\mathbf{x}^s_{(j)}(\mathbf{x}_i)$ as the point which is the $j$th point closest to the sphere of the $d$-dimensional ball with center at $\mathbf{x}_i$ and radius equal to $s$, i.e.,

$$\left| s - ||\mathbf{x}^s_{(1)}(\mathbf{x}_i) - \mathbf{x}_i|| \right| \leq \cdots \leq \left| s - ||\mathbf{x}^s_{(j)}(\mathbf{x}_i) - \mathbf{x}_i|| \right|$$
$$\leq \cdots \leq \left| s - ||\mathbf{x}^s_{(n-1)}(\mathbf{x}_i) - \mathbf{x}_i|| \right|. \qquad (1)$$

Let $\mathbf{p}^1_i(s)$ be the line which passes through the two points $\mathbf{x}_i$ and $\mathbf{x}^s_{(1)}(\mathbf{x}_i)$, i.e., $\mathbf{p}^1_i(s) = \left\{ \mathbf{y} : \mathbf{y} - \mathbf{x}_i = t_1(\mathbf{x}^s_{(1)}(\mathbf{x}_i) - \mathbf{x}_i), \forall t_1 \in \mathrm{R} \right\}$. Let $\theta^1_i(s)$ be the angle between the vector $\mathbf{x}^s_{(2)}(\mathbf{x}_i) - \mathbf{x}_i$ and its orthogonal projection on $\mathbf{p}^1_i(s)$. Hence $\theta^1_i(s)$ only take values between 0 and $\frac{\pi}{2}$.

Based on the collection of $\theta^1_i(s)$, $i = 1, \ldots, n$ over the whole data set, the following heuristic is useful for any 1-d manifold. Since a 1-d manifold is a topological space which is locally Euclidean, then $\theta^1_i(s)$ should be close to 0 for small values of $s$ if the data lie on a 1-d manifold. When the data have noise, the degree of tolerance of the noise for all window sizes depends on the different noise levels. For the first toy example with small noise, both scales are big enough to perceive the 1-d structure since the averages of $\theta^1_i(s)$ at the small and big scales are both close to 0 ($\overline{\theta^1(s)}$ equal to 0.013 and 0.004 respectively). For the second toy example, the value of $\overline{\theta^1(s)} = 0.599$ is much larger at the small scale, reflecting the virtually apparent 2-d structure, and $\overline{\theta^1(s)} = 0.112$ is smaller for larger $s$, reflecting the less 2-d nature shown in Figure 2. $\overline{\theta^1(s)}$ of the second toy example indicates more noise than the first one. For the third toy example, both scales suggest 2-d manifold structure because their averages are far from 0 (0.724 and 0.738).

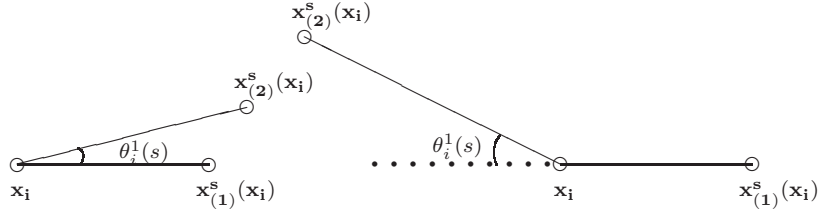

FIG 4. *The two cases of the definition of the angle $\theta_i^1(s)$ between $\mathbf{x}_{(1)}^s(\mathbf{x_i})-\mathbf{x_i}$ and $\mathbf{x}_{(2)}^s(\mathbf{x_i})-\mathbf{x_i}$, always taken to be between 0 and $\frac{\pi}{2}$.*

Theoretically, we could have $\theta_i^1(s)$ for any scale $s$. Notice that for a given data set, $\theta_i^1(s)$ for $i = 1,\ldots,n$ are the same when $s \leq s_{\min}$ and $s \geq s_{\max}$, where $s_{\min}$ and $s_{\max}$ are the minimum and the maximum pairwise distances among the data points. Indeed, when $s \leq s_{\min}$ is considered, for each data point $\mathbf{x}_i$, $\theta_i^1(s)$ will not be changed since the 2 nearest points are always the same. Similarly, when $s \geq s_{\max}$, $\theta_i^1(s)$ will not be changed for the same reason. Therefore, for a given data set, we only need to concentrate on the scales $s$ between $s_{\min}$ and $s_{\max}$. Divided by $s_{\max}$, we will have the standardized scale from 0 to 1.

We would like to explore the "tendency of a data set to lie on or nearly to lie on a 1-d submanifold" by analyzing the set of angles, $\theta_i^1(s)$. However, these quantities are not very interpretable. A sensible summary of $\theta_i^1(s)$ is the average of $\theta_i^1(s)$. A more interpretable statistic $T_1(s)$ can be defined as a re-scaling of $\overline{\theta^1(s)}$,

$$T_1(s) = 1 + \frac{1}{a_1} \times \overline{\theta^1(s)} \qquad (2)$$

where $a_1 = \frac{\pi}{4}$. Such an $a_1$ value is chosen to match $T_1(s)$ as the dimension of the data. Details are explained later in this section.

**Theorem 2.1.** *Assume that $\mathbf{Y} = \{X_1,\ldots,X_n\}$ is a random sample and the density function is absolutely continuous with respect to Lebesgue measure over an appropriate manifold. Then the following are true.*

- *if $\mathbf{Y}$ is from a 1-d linear space, $P(\theta_i^1(s) \equiv 0) = 1$, from a 1-d manifold, $P(\theta_i^1(s_{min}) = 0) \to 1$ as $n \to \infty$;*
- *if $\mathbf{Y}$ is from a 2-d manifold, $\theta_i^1(s_{min}) \xrightarrow{L}$ (the uniform$(0, \frac{\pi}{2})$ distribution) as $n \to \infty$. Further $E(\theta_i^1(s_{min})) \to \frac{\pi}{4}$ as $n \to \infty$;*
- *if $\mathbf{Y}$ is from a $k(k > 2)$-d manifold, $\lim_{n\to\infty} E(\theta_i^1(s_{min}))$ increases as $k$ increases, with limit equal to $\frac{\pi}{2}$ as $k \to \infty$.*

Therefore, $T_1(s)$ can be used to distinguish a data set from a manifold with dimensionality $\leq 2$ or $> 2$. This idea can be generalized to distinguish a data set with dimensionality $\leq k+1$ or $> k+1$, where $k > 1$. The idea of angles can be extended by considering $(k+2)$-tuples. For each point $\mathbf{x_i}$ and a fixed value $s$, consider the $k+1$ points which are closest to the sphere of the $d$-dimensional ball with center at $\mathbf{x}_i$ and radius equal to $s$, denoted as $\mathbf{x}^s_{(1)}(\mathbf{x}_i), \ldots, \mathbf{x}^s_{(k+1)}(\mathbf{x}_i)$ as in (1). Let $\mathbf{p}^k_i(s)$ be the hyperplane determined by the points $\mathbf{x}_i, \mathbf{x}^s_{(1)}(\mathbf{x}_i), \ldots, \mathbf{x}^s_{(k)}(\mathbf{x}_i)$. Then the angle $\theta^k_i(s)$ is formed by the vector $\mathbf{x}^s_{(k+1)}(\mathbf{x}_i) - \mathbf{x}_i$ and its orthogonal projection on $\mathbf{p}^k_i(s)$. Therefore, for a data set with $n$ points, from a $d$-dimensional space, we can consider effective dimensions $k = 1, \ldots, m$ based on the angles: $\mathbf{\Theta}(\mathbf{s}) = \{\theta^j_i(s), i = 1, \ldots, n; j = 1, \ldots, m\}$, where $m = \min(d-1, n-2)$. Again all of these angles are defined to be between $0$ and $\frac{\pi}{2}$.

Each set $\boldsymbol{\theta}^k(s) = \{\theta^k_i(s), i = 1, \ldots, n\}$ will help to determine whether the data are close to a manifold with dimensionality $\leq k+1$ or $> k+1$. Correspondingly, we can define the more interpretable statistics $\mathbf{T}(\mathbf{s}) = \{T_k(s), k = 1, \ldots, m\}$ as follows

$$T_k(s) = k + \frac{1}{a_k} \overline{\theta^k(s)}, \qquad (3)$$

where for a fixed point $\mathbf{x}_i$, $a_k$ is the mean of the limit distribution of $\theta^k_i(s_{\min})$ as $n \to \infty$ for the data from the $(k+1)$-d manifold. Denote $\eta_k$ be the random variable with such a limit distribution, so $a_k = \mathrm{E}\left(\eta^k\right)$. For example, $a_1 = \frac{\pi}{4}$. The limit distributions of all $\theta^k_i(s_{\min})$ $(k > 1)$ are given in Chapter 5 of [27].

## 3. Effective dimensionality

In the previous section, we studied the statistics $\mathbf{\Theta}(s)$ and $\mathbf{T}(s)$ which are both based on the data. In this section, we are going to define the definition for a population $\mathcal{X}$ analogue to $\mathbf{T}(s)$, denoted as $\mathbf{D}(s)$. We call this $\mathbf{D}(s)$ the *Vector Dimensionality* of the population. Under some general conditions, we will show the consistency of the statistic $\mathbf{T}(s)$ for $\mathbf{D}(s)$. Finally, we will set up hypothesis tests to extract the effective dimensionality of data based on $\mathbf{T}(s)$.

### 3.1. Definition of Vector Dimensionality

Assume that a population, $\mathcal{X}$, is from a $d$-dimensional space, $\mathcal{R}^d$, and has a probability distribution $\mathbf{P}(x)$. For a random point $\mathbf{x} \in \mathcal{X}$, if we use polar coordinates and set the origin of the coordinate system at $\mathbf{x}$, then for every $\mathbf{y} \in \mathcal{X}$, we will have polar coordinates $(\rho, \xi_1, \ldots, \xi_{d-1})$, simply denoted as $(\rho, \Lambda)$ where $\Lambda = (\xi_1, \ldots, \xi_{d-1})$. For any scale $s \in (0, \infty)$, let $\mathcal{Y}|_{(\mathbf{x},s)}$ be the conditional population of points $\mathbf{y}$ on the surface of the ball $S_{\mathbf{x},s}$, i.e., $\mathbf{y} \in \mathcal{Y}|_{(\mathbf{x},s)}$ satisfies $||\mathbf{y} - \mathbf{x}|| = s$. For every $\mathbf{y} \in \mathcal{Y}|_{(\mathbf{x},s)}$, $\mathbf{y}$ will have polar coordinates $(s, \Lambda)$. And the distribution of points in the population $\mathcal{Y}|_{(\mathbf{x},s)}$ will be a conditional probability measure of $\Lambda$ given $\rho = s$ and the point $x$, denoted as $\mathbf{P}_x(\cdot|\rho = s)$.

For any two random points $\mathbf{y}_1, \mathbf{y}_2 \in \mathcal{Y}|_{(\mathbf{x},s)}$, let $\mathbf{L}_{\mathbf{x}}^1$ be the line determined by $\mathbf{x}$ and $\mathbf{y}_1$, then $\theta_{\mathbf{x}}^1(s)$ will be the angle between the vector $\mathbf{y}_2 - \mathbf{x}$ and its projection on $\mathbf{L}_{\mathbf{x}}^1$. In fact, $\theta_{\mathbf{x}}^1(s)$ is a function only of $\Lambda_1$ and $\Lambda_2$, denoted as $h_1$, i.e., $\theta_{\mathbf{x}}^1(s) = h_1(\Lambda_1, \Lambda_2)$. Define the moments:

$$\psi_{\mathbf{x}}^1(s) = \int h_1(\Lambda_1, \Lambda_2) \, d\mathbf{P}_x(\Lambda_1|\rho = s) d\mathbf{P}_x(\Lambda_2|\rho = s).$$

Further more, in general we can define the moments:

$$\psi_{\mathbf{x}}^{k-1}(s) = \int h_{k-1}(\Lambda_1, \ldots, \Lambda_k) \, d\mathbf{P}_x(\Lambda_1|\rho = s) \cdots d\mathbf{P}_x(\Lambda_k|\rho = s), \quad \text{for } k < d.$$

In order to get properties of the whole population, we need to consider every point, i.e., we need to consider $\int \psi_{\mathbf{x}}^1(s) d\mathbf{P}(\mathbf{x})$, ..., $\int \psi_{\mathbf{x}}^{d-1}(s) d\mathbf{P}(\mathbf{x})$. To make them more interpretable, we will re-scale $\int \psi_{\mathbf{x}}^1(s) d\mathbf{P}(\mathbf{x})$, ..., $\int \psi_{\mathbf{x}}^{d-1}(s) d\mathbf{P}(\mathbf{x})$ to a series $D_1(s)$, ..., $D_{d-1}(s)$ so that: when the population is contained in a $d'$-dimensinal linear space, we have $D_1(s)$, ..., $D_{d-1}(s)$ as in Table 1.

TABLE 1
Values of $D_k(s)$ for linear spaces

| $d'$ | $D_1(s)$ | $D_2(s)$ | $\cdots$ | $D_k(s)$ | $\cdots$ | $D_{d-2}(s)$ | $D_{d-1}(s)$ |
|---|---|---|---|---|---|---|---|
| 1 | 1 | 2 | $\cdots$ | $k$ | $\cdots$ | $d-2$ | $d-1$ |
| 2 | $\geq 1$ | 2 | $\cdots$ | $k$ | $\cdots$ | $d-2$ | $d-1$ |
| $k$ | $\geq 1$ | $\geq 2$ | $\cdots$ | $k$ | $\cdots$ | $d-2$ | $d-1$ |
| $d-1$ | $\geq 1$ | $\geq 2$ | $\cdots$ | $\geq k$ | $\cdots$ | $\geq d-2$ | $d-1$ |

Now we have the following definition of the population version of the sample statistic $\mathbf{T}(s)$.

**Definition 3.1.** Vector Dimensionality $\mathbf{D}(s)$ of a population $\mathcal{X}$ in $\mathcal{R}^d$ with a probability $\mathbf{P}$ is defined as $\mathbf{D}(s) = \{ D_1(s), D_2(s), \ldots, D_{d-1}(s) \}$, where for any fixed value $s \in (0, \infty)$

$$D_1(s) = 1 + \frac{1}{a_1} \int \psi_{\mathbf{x}}^1(s) d\mathbf{P}(\mathbf{x}), \tag{4}$$

$$\ldots$$

$$D_{d-1}(s) = d - 1 + \frac{1}{a_{d-1}} \int \psi_{\mathbf{x}}^{d-1}(s) d\mathbf{P}(\mathbf{x}). \tag{5}$$

where $a_1$, ..., $a_{d-1}$ are defined as before.

### 3.2. Consistency of the statistic $\mathbf{T}(s)$

The natural statistic for $\mathbf{D}(s)$ is $\widehat{\mathbf{D}} = \mathbf{T}(s) = \{ T_1(s), T_2(s), \ldots, T_m(s) \}$, where $m = \min(d - 1, n - 2)$. In this section, we are going to prove that for any $k = 1, \ldots m$, $T_k(s)$ is a consistent estimate for $D_k(s)$ for any fixed value of $s$.

The following lemma is useful in proving the consistency of $\mathbf{T}(s)$. The idea of the proof is similar to the proof of Lemma 5.1 in Section 5.2 in Devroye, Györfi, and Lugosi (1996).

Assume $X_1, X_2 \ldots$ is a random sequence from a population $\mathcal{X}$ with a probability measure P on it. Let $S_{x,\epsilon,s}$ be the set $\{\, y \in S_{x,\epsilon,s} : s-\epsilon \leq \|y-x\| \leq s+\epsilon \,\}$, for any $\epsilon < s$. Because the probability measure P has a continuous density function $f$, it follows that $\mathrm{P}\,(S_{x,\epsilon,s}) > 0$ for any $x \in \mathcal{X}$, $\epsilon > 0$, and $s > 0$. For $X_1, X_2, \ldots, X_n$, order $|\|X_i - x\| - s|$ and define analogues of order statistics $X_{(k)}^s(x)$ such that

$$\left|\|X_{(1)}^s(x) - x\| - s\right| \leq \left|\|X_{(2)}^s(x) - x\| - s\right| \leq \cdots \leq \left|\|X_{(n)}^s(x) - x\| - s\right|.$$

So $X_{(k)}^s(x)$ is the point among $X_1, X_2 \ldots, X_n$, which is the $k$-th closest to the surface of the closed ball centered at $x$ with radius $s(> 0)$. Then we have the following lemma.

**Lemma 3.1.** *For any $x \in \mathcal{X}$ and any fixed value of $k$, $\left|\|X_{(k)}^s(x) - x\| - s\right| \to 0$ with probability 1 as $n \to \infty$. If $X$ is independent of the data, then*

$$\left|\|X_{(k)}^s(X) - X\| - s\right| \to 0$$

*with probability 1 whenever $n \to \infty$.*

**Theorem 3.2.** *Assume a random sample $\mathbf{Y} = \{\, X, X_1, \ldots, X_n \,\}$ from the population $\mathcal{X} \subset \mathcal{R}^d$, and $\theta_X^k(s)$ defined as above, then $E\theta_X^k(s)$ converges to $\psi_X^k(s)$ as $n \to \infty$. Further more, for any fixed value $s$ and any fixed integer $k$, $\lim_{n \to \infty} ET_k(s) = D_k(s)$.*

**Lemma 3.3.** *Suppose the observations $\{\, X_i, i = 1, \ldots, n \,\}$ are independent, identically distributed in d-dimensional space with a twice continuously differentiable density $f$. Define $T_k(s)$ as (3), then for any fixed value of scale $s$ and any fixed integer $k$, where $1 \leq k \leq m$, $Var(T_k(s)) \to 0$ as $n \to \infty$.*

**Theorem 3.4.** *For any fixed value of scale $s$ and fixed integer $k$, where $1 \leq k \leq m$,*

$$T_k(s) \xrightarrow{P} D_k(s), \quad as\ n \to \infty$$

The proof of the above lemmas and theorems in this section are given in [27].

### 3.3. Hypothesis test

$\mathbf{D}(s)$ is a theoretical construction that allows understanding of our approach and estimating effective dimensionality. However, we generally do not know the distribution of the true population. Instead, we use the sample statistic $\mathbf{T}(s)$. A set of hypothesis tests is used, each one based on one element of $\mathbf{T}(s)$. Each test is performed separately. Here, the hypothesis test is a mechanism to extract the effective dimensionality of the data. A simple approach to the multiple comparison issue would be a Bonferroni method. More sophisticated multiple comparison tests for effective dimensionality are interesting topics for the future work.

For the typical $k$-d data, we will assume the $k$-dimensional standard normal distribution as the distribution under the null hypothesis. Since analytical calculation of $D_k(s)$ of the $(k+1)$-d standard normal distribution seems to be intractable, we use a simulation approach. For example, for $D_1(s)$ of the 2-d standard normal distribution, denoted as $D_1^2(s)$, we generated 1000 random samples with each sample of size equal to the sample size of the testing data, and calculated $T_1(s)$ for these 1000 random samples. We use the average of these 1000 $T_1(s)$ as an estimate of $D_1^2(s)$, denoted as $\widehat{D_1^2(s)}$, and use the 2.5% and 97.5% percentiles of these 1000 $T_1(s)$, denoted as $\widehat{D_{1,2.5\%}^2}(s)$ and $\widehat{D_{1,97.5\%}^2}(s)$, as the critical value at the 5% significance level. Since $D_1(s) \equiv 1$ for any 1-d linear data, by comparing whether $T_1(s)$ of the data less than $\widehat{D_{1,2.5\%}^2}(s)$ or not, we will conclude for the fixed scale value $s$, whether the effective dimension of the data is less than 2-d at the 2.5% significant level or not. If $T_1(s)$ of the data is greater than $\widehat{D_{1,2.5\%}^2}(s)$ and less than $\widehat{D_{1,97.5\%}^2}(s)$, then the effective dimension for that scale $s$ is 2 at the significant level of 5%. The steps, to investigate the effective dimensionality of data at all scales, follow. For each $s \in S$, or in particular a discrete subset of $S$,

1. For $k = 1, 2, \ldots$, test

    H$_0$  $T_k(s)$ of the data is significantly equal or less than $\widehat{D_k^{k+1}}(s)$ for $s \in S$ at level $\alpha$

    vs. H$_1$  otherwise
2. stop when H$_0$ is accepted for this scale $s$, and record the corresponding $T_k(s)$ of the data as the effective dimension of the testing data set at this scale $s \in S$.

## 4. Examples

In this section we present several examples with applications of our method. The first two examples are based on simulated data. Example 3 and 4 are two real data sets. Example 5 is about the deterministic chaos case.

### *4.1. Simulated examples*

**Example 1: A circle**

We use Example 1 and 2 in "Adaptive Principal Surfaces" by LeBlanc & tibshirani, (1994). Generate the first data set of 100 observations from the circle

$$y_1 = 5\sin(\lambda) + \epsilon_1$$
$$y_2 = 5\cos(\lambda) + \epsilon_2,$$

where $\lambda \sim \mathrm{U}(0, 2\pi)$ and $\epsilon_i \sim \mathrm{N}(0, 0.25)$ for $i = 1, 2$. And by adding 4-dimensional noise, $y_i = \epsilon_i$, where $\epsilon_i$, $i = 3, 4, 5, 6$ have the same distribution as $\epsilon_1$, we have the second data set.

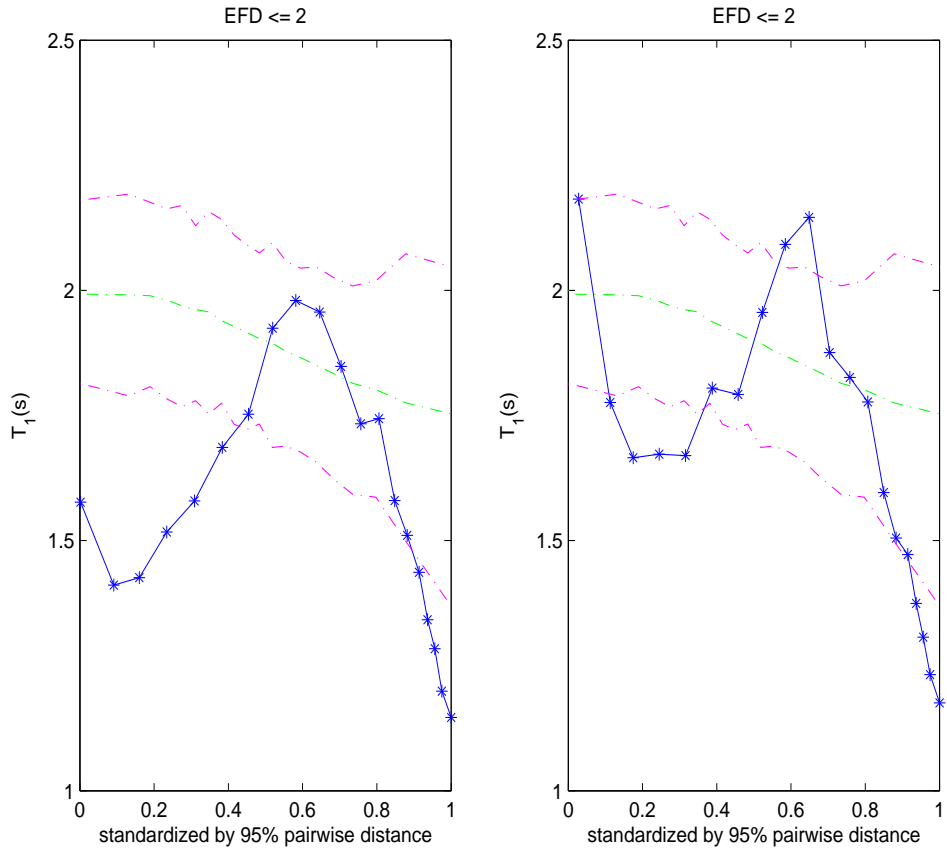

Fig 5. *Plot of $T_1(s)$ for the circle examples: 2-d (left) and 6-d (right)*

In plots of Figure 5, we compare $T_1(s)$ based on the two simulated data sets with the $T_1(s)$ of the 2-d standard normal distribution. The green dash-dotted curve and the 2 magenta dash-dotted curves represent the estimate of $T_1(s)$ and the "typical range" (95% confidence bands) of $T_1(s)$ if the data follow a 2-d standard normal distribution. The blue "-*" curves are $T_1(s)$ based on the simulated data sets. The scale $s$ is chosen from the minimum pairwise distance to the maximum increased by 5% of pairwise distance, denoted as $s_{\min}, s_{5\%}, s_{10\%}, \ldots, s_{95\%}, s_{\max}$. To make the scale comparable, we standardize it either by the maximum pairwise distance or by the 95% pairwise distance to reduce the affect of outliers.

For the first simulated data set, $T_1(s)$ is less than lower magenta line for relatively small and large scales, indicating that the effective dimension is leas than 2 for those scales. For some middle range scales, $T_1(s)$ is within the 2 magenta lines, indicating that they are 2-d. $T_1(s_{\min}) > T_1(s_{5\%})$ because the data have noise and the noise affect the effective dimension more at the smaller scales. $T_1(s)$ reaches its maximum value at some middle range of scale for the reason that the underlying structure of the data is nonlinear and shows 2-d structure at the middle scales. For the second data set, $T_1(s)$ is larger than the first case at every scale because of the larger noise level. However, it still shows less than 2-d structure of the data at some scales, see right plot of Figure 5. $T_1(s)$ is significant bigger than 2 when $s = s_{45\%}$, $s_{50\%}$. We should go further to test $T_2(s)$ for this data set at these two scales to determine the effective dimensionality at these two scales. By applying ISOMAP to these two data sets, no matter what the value of neighborhood parameter is taken, ISOMAP always estimates the intrinsic dimension as 2 or bigger.

**Example 2: The Swiss roll**

In this example, we will compare our method with ISOMAP through two toy data sets, the noiseless Swiss roll data in plot (a) of Figure 1 and the noisy Swiss roll data in plot (b) of Figure 1.

For the noiseless Swiss roll data shown in plot (a) of Figure 1, the result analyzed by the ISOMAP algorithm is shown in plot (a) of Figure 6. As with PCA, the true dimensionality of the data can be estimated from the decrease in residual variance as the dimensionality of the low-dimensional representation increases. The residual variance of ISOMAP correctly bottoms out at dimensionality equal to 2.

Since the data lie exactly on the surface of the Swiss roll, the shapes of $T_1(s)$ and $T_2(s)$ are driven completely by the curvature of the surface of the Swiss roll. Plot (b) in Figure 6 shows $T_1(s)$ of the noiseless data (blue "-*" line) and compares it with the $T_1(s)$ of the typical 2-d data (the typical data means the data with standard normal distribution). $T_1(s)$ is significantly greater than the typical $T_1(s)$ of 2-d data from $s_{\min}$ to $s_{75\%}$. It suggests that the data are close to a manifold with dimensionality greater than 2. From $T_2(s)$ in plot (c), at $s_{\min}$, $s_{5\%}$, and $s_{10\%}$, the corresponding values of $T_2(s)$ are significantly less than the typical $T_2(s)$ of 3-d data, especially $T_2(s_{\min})$ is about 2.2 and $T_2(s_{\min}) < T_2(s_{5\%})$. Since the noise affects detection of the real structure mostly at the smallest scale, the fact that $T_2(s_{\min}) < T_2(s_{5\%})$ indicates that even at the smallest scale, the window sizes is big enough to detect the true underline structure. $T_2(s)$ of this simulated data set starts significantly less than the typical $T_2(s)$ of 3-d data at $s_{\min} \leq s < s_{15\%}$, increases to values bigger than the typical $T_2(s)$ of 3-d data at the medium scales, then decreases to the values less than the typical $T_2(s)$ of 3-d data again at large scales. The values of $T_2(s)$ consistently less than 3 at smaller scales indicate that the data are close to 2-d manifold. The non-monotone shape suggests that the data lie on a curved 2-d manifold, since the 3-d structure only shows at the medium scales.

For the noisy Swiss roll data in plot (b) of Figure 1, plot (d) of Figure 6 shows the estimated dimensionality by ISOMAP. We tried the values of neighborhood size from 2 to 22 by every increment of 2. The result shown in plot (d) is the best among them. Plots (e) and (f) show our test statistics $T_1(s)$ and $T_2(s)$ of the noisy Swiss roll data comparing to the typical 2-d and 3-d data. Our method suggests that this noisy Swiss roll data are still close to a 2-d manifold but with relatively large noise, since $T_2(s_{\min})$ and $T_2(s_{5\%})$ are significantly less than the corresponding typical $T_2(s)$ of 3-d data. However, comparing to the result of noiseless Swiss roll data, $T_2(s_{\min}) \approx 2.85$ is much bigger than the previous one, $T_2(s_{\min}) \approx 2.2$.

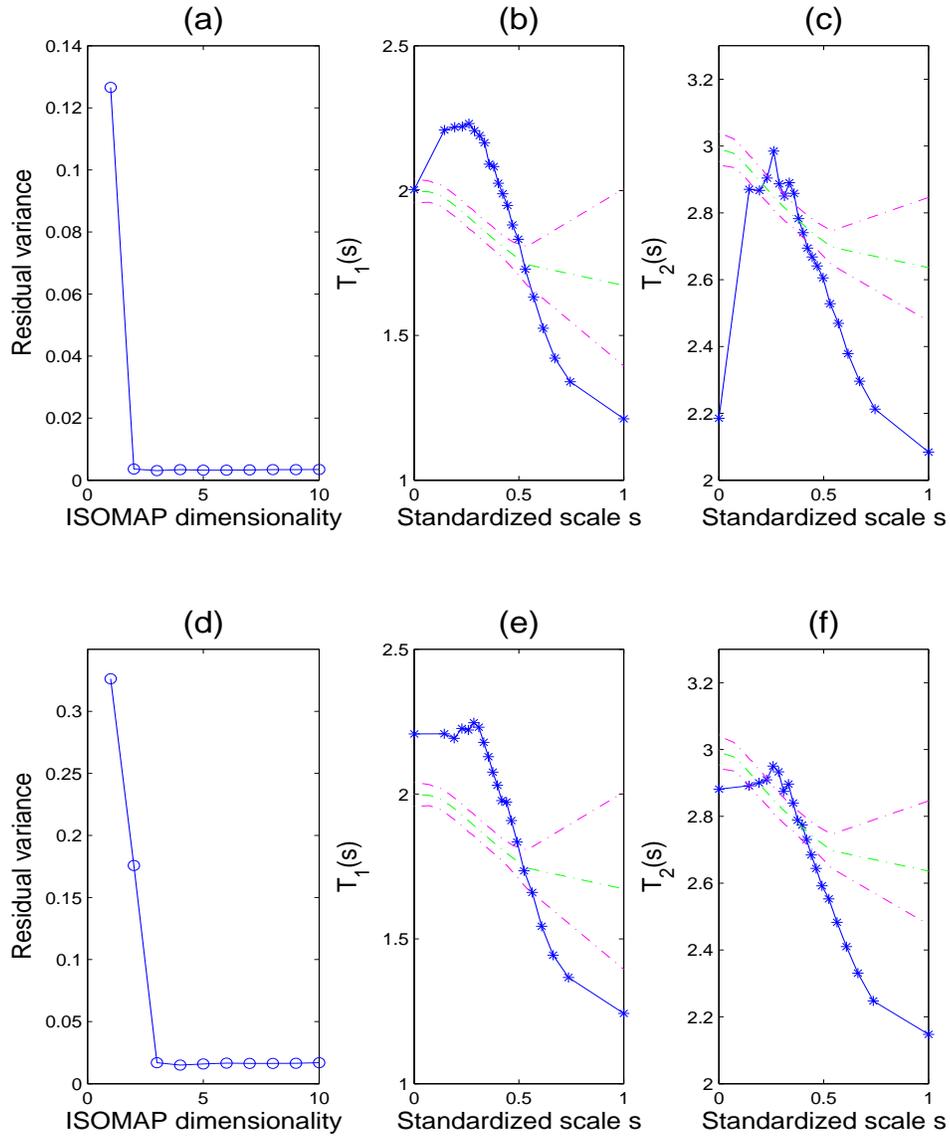

FIG 6. *(a): The estimated dimension by ISOMAP for the noiseless Swiss roll data in Figure 1. (b): For the noiseless Swiss roll data: comparing their $T_1(s)$ with the typical $T_1(s)$ of 2-d data. (c): For the noiseless Swiss roll data: comparing their $T_2(s)$ with the typical $T_2(s)$ of 3-d data. (d): The estimated dimension by ISOMAP for the noisy Swiss roll data in Figure 1. (e): For the noisy Swiss roll data: comparing their $T_1(s)$ with the typical $T_1(s)$ of 2-d data. (f): For the noisy Swiss roll data, comparing their $T_2(s)$ with the typical $T_2(s)$ of 3-d data.*

### *4.2. Two real data sets*

**Example 3: Berkeley growth study data**

The Berkeley growth study data (Chapter 6, Ramsay and Silverman 2002) have been studied by Wang and Iyer (2006). The data set contains height measurements of 54 girls and 39 boys with total 31 measurements on each subject from age 1 to 18 years. For comparison reason, we choose the last 21 measurements from age 8 to 21 years which were analyzed by Wang and Iyer, i.e., a data set with 103 observations on 21-dimensional space. Wang and Iyer developed a method of finding nonlinear latent structure by LLE. To find the intrinsic dimension, they use the cross-validation method. The left panel of Figure 7 shows that such a data set has a significant 2-d underlying structure which is consistent to the result of Wang and Iyer. They concluded that 2 factors essentially determine the growth patten for the juvenile boys and girls.

**Example 4: States data**

The data consist of 7 variables for 50 U.S. states, which are population, average income, illiteracy rate, life expectancy, homicide rate, high school graduation rate, and average number of days with below-freezing minimum temperatures. The data are available from Becker, Chambers, and Wilks (1988) and were used as an example in the paper of LeBlanc and Tibshirani (1997). In their paper, they showed that the projection onto the 2-dimensional adaptive principal surface captured the most variation of the data. As we explained before, they use GCV to find the intrinsic dimension. The right panel of Figure 7 shows that for the most of scale, the data are lying closely on the 2-dimensional manifold.

From Example 3 and 4, it is clear that we can find the effective dimensions of these two data sets very quickly based on out scale space approach.

### *4.3. Intrinsic dimension for deterministic chaos*

Deterministic chaos has been rigorously and extensively studied by mathematicians and scientists. Instead of presenting a formal account, we will adopt an informal approach in which we illustrate some basic concepts of deterministic chaos through an example.

**Example 5: Hénon map**

The example is the Hénon map defined by

$$\begin{aligned} x_n &= y_{n-1} + 1 - 1.4x_{n-1}^2 \\ y_n &= 0.3x_{n-1} \end{aligned} \quad (6)$$

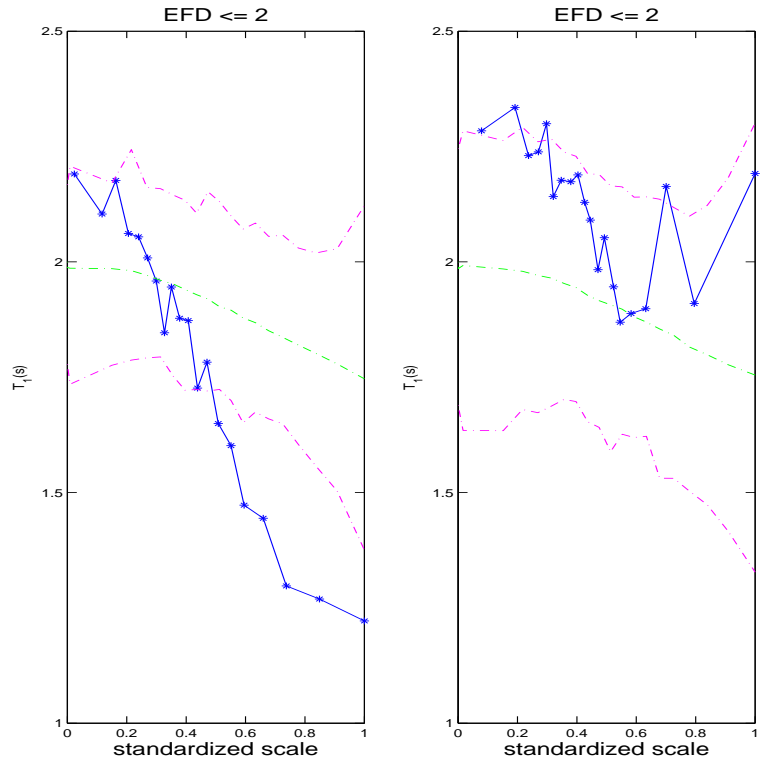

FIG 7. *The left panel is the plot of $T_1(s)$ for the Berkeley growth study data. The right panel is the plot of $T_1(s)$ for the States data.*

which maps a 2-d vector $(x_{n-1}, y_{n-1})$ to $(x_n, y_n)$. There is an extensive literature about the Hénon map. However, it has proved surprisingly difficult to obtain rigorous results for this map. The correlation dimension, $p$, has been quoted as being in the range $1.25 \pm 0.02$ (Grassberger and Procaccia, 1983). The direct Grassberger-Procaccia method, the Takens's estimate, and other modified Grassberger-Procaccia methods (they are essentially equivalent) give similar estimated values $\hat{p} = 1.22$. When noise is present, all these methods are not effective.

In Smith's (1992a) paper, there is a theoretical discussion about the presence of noise and a modified version of the Grassberger-Procaccia methods for estimating the dimensionality in the presence of the noise. In his paper, he used the Hénon map as an example. Here we follow Smith's example and generate 1100 data points from the Hénon map as (6) by discarding the first 100 and adding the 2-d Gaussian observational noise with $\sigma = 0, 0.001, 0.003, 0.01$. Figure 8 shows the four data sets. The $\sigma = 0.001$ data set appears no different from the case $\sigma = 0$, whereas at $\sigma = 0.003$ the difference is noticeable but the fractal structure still evident, and at $\sigma = 0.01$ the picture has become quite seriously blurred.

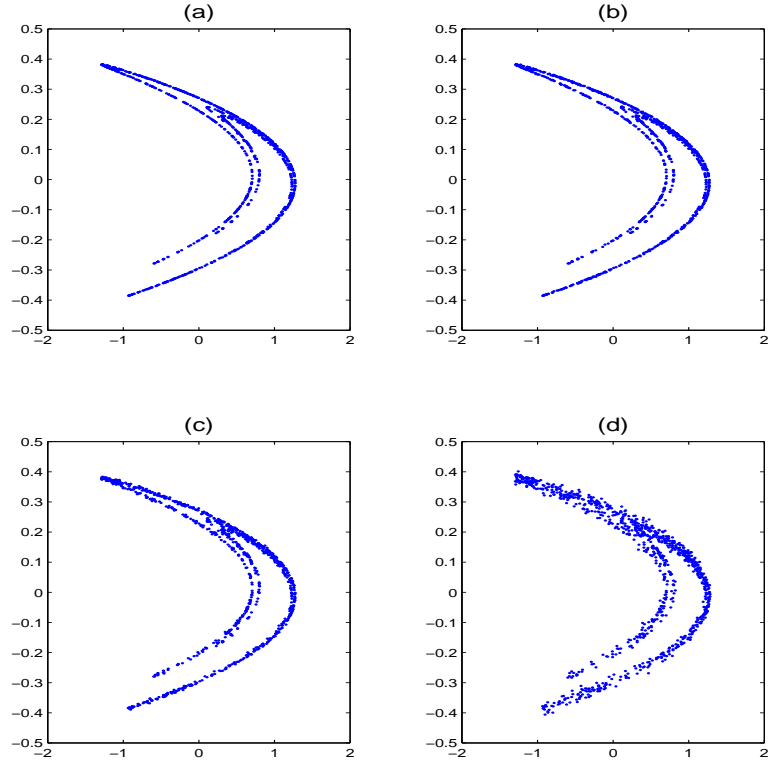

FIG 8. *The scatter plots for the trajectory starting from the origin and driven by the Hénon map plus the 2-d Gaussian noise. Sample size is 1000. (a) with no noise, i.e., $\sigma = 0$; (b) with $\sigma = 0.001$; (c) with $\sigma = 0.003$; (d) with $\sigma = 0.01$.*

For the data set with $\sigma = 0.001$, recall that this is visually indistinguishable, some estimates of $p$ are as high as 2.5. For the data set with $\sigma = 0.003$, the estimates of $p$ vary from 1.3 to 2.7. For the data with $\sigma = 0.01$, the estimates are even worse. Details can be found in Smith (1992b). In that paper, Smith gave a dimension estimate dealing with the noise, which requires either the prior knowledge or an estimation of the variance of the noise. Smith's method gave much improved estimates of the dimensions of the data with noise. Because we study the case without prior knowledge, we will only compare our method with the methods which require no knowledge of or no estimation of the variance, *e.g.,* Takens' estimate.

Figure 9 shows $T_1(s)$ for the four simulated data sets. At $s_{\min}$, $T_1(s)$ increases as the noise level increases. This is not surprising, since the smallest scale senses the high-dimensional structure. However, $T_1(s_{5\%})$ is around 1.4 for the data sets with $\sigma = 0$, $\sigma = 0.001$, and $\sigma = 0.003$. $T_1(s_{5\%})$ is about 1.5 for the data set with $\sigma = 0.01$. Our method shows that all 4 data sets are significantly different from a 2-d manifold.

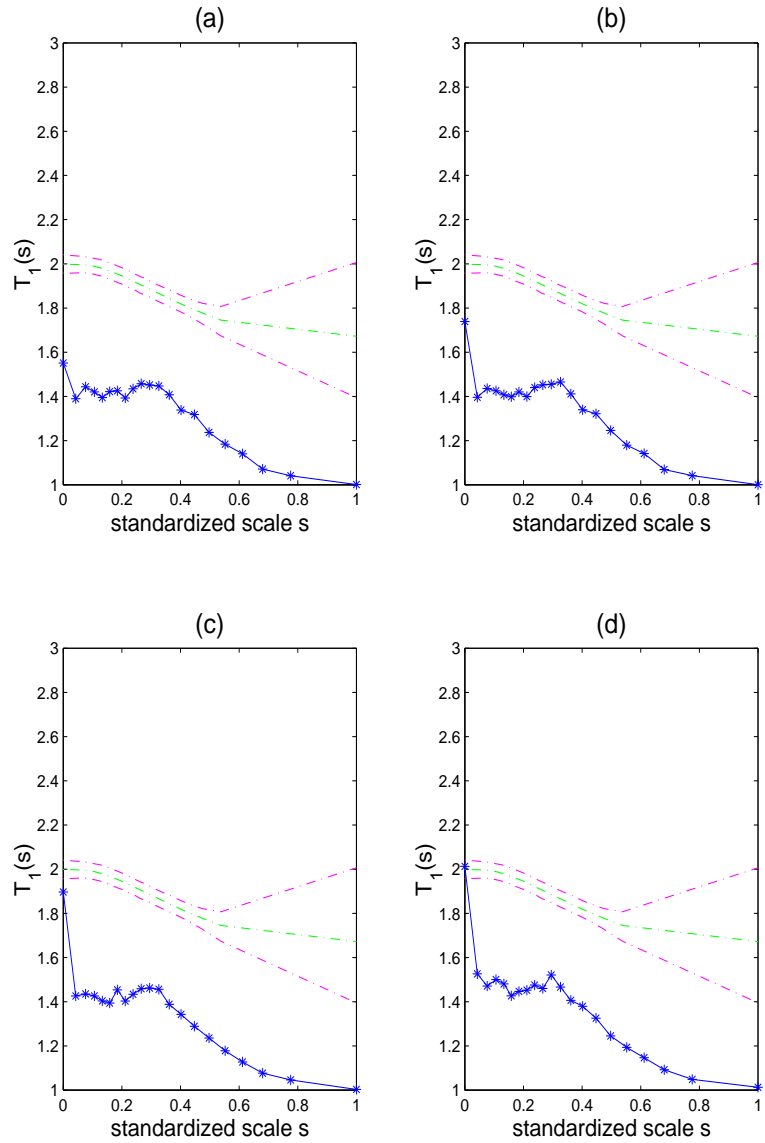

Fig 9. *Comparison between $T_1(s)$ of the 4 data sets in Figure 8 and $T_1(s)$ of the typical 2-d data. They are all significantly different from a 2-d manifold.*

To summarize the comparison between our method and the methods for estimating the correlation dimension of the deterministic chaos, the similarity is that both two versions of estimating dimensionality try to find the effective dimensionality between integers. In other words, they not only give the information whether the data are close to a $k$-dimensional manifold with $k$ as an integer, but also provide the information about the degree of the closeness by using fractional dimensionality. They give similar information when the data are studied locally. However, our method gives the effective dimensionality not only locally but also globally. Unlike the general methods of estimating the correlation dimension for the deterministic system, the noise in the data does not affect the stability of our method.

## 5. Summary and discussion

The method provided in this paper does not give a single number as an estimate of intrinsic dimension as others. In stead, it gives a whole picture of intrinsic dimension for data at all scales. The estimated intrinsic dimension is a function of scale parameter $s$ which takes values not integer but real number. The fraction part of the intrinsic dimension of the small scales helps us to learn some information about the noise. Based on our method, it is possible to choose a proper estimate for the intrinsic dimension as an input for many nonlinear dimension reduction methods. Our procedure does not depend on any specific dimension reduction method. It is more robust than ISOMAP when data have relatively large noise and more computational efficient than most of CV and GCV methods.

A couple of extensions of the procedure might be worth of investigating.

- For a fixed scale $s$ and a fixed point, only two points at this scale are taken into account to calculate the test statistics. It could be extended to include more data points at the scale so that the test statistics are more robust.
- Currently, we use the simulation result to determine the significance of the effective dimension. It is important to establish the asymptotic distributions of test statistics at all scale.

Our current method is carried out in a simple Matlab code.